\documentclass[12pt]{article}%
\usepackage{amsmath}
\usepackage{amsfonts}
\usepackage{amssymb}
\usepackage{graphicx}%
\providecommand{\U}[1]{\protect\rule{.1in}{.1in}}
\begin{document}

\title{On Noncommutative Levi-Civita Connections}
\date{}
\author{Mira A. Peterka\thanks{Partially supported by the Simons Foundation grant
346300 and the Polish Government MNiSW 2015-2019 matching fund.}
\thinspace\ and Albert J. L. Sheu}
\maketitle

\begin{abstract}
We make some observations about Rosenberg's Levi-Civita connections on
noncommutative tori, noting the non-uniqueness of general torsion-free
metric-compatible connections without prescribed connection operator for the
inner *-derivations, the nontrivial curvature form of the inner *-derivations,
and the validity of the Gauss-Bonnet theorem for two classes of non-conformal
deformations of the flat metric on the noncommutative two-tori, including the
case of non-commuting scalings along the principal directions of a two-torus.

\end{abstract}

\begin{description}
\item[Keywords: ] Noncommutative tori, Levi-Civita connection, Riemannian
curvature, Gauss-Bonnet theorem.

\item[Mathematics Subject Classification (2010): ] 46L87, 58B34,
46L08.\newline\newline
\end{description}

\section{Introduction}

In the recent paper \cite{Ro}, Jonathan Rosenberg revisits the concept of
connections on noncommutative tori \cite{Ri2} introduced by Connes
\cite{C1,C2}, a concept that resulted in the classification of Yang-Mills
moduli for the noncommutative two-torus \cite{CRi,Ri1}, and also features
prominently in the vast literature of noncommutative gauge theory (see for
example the references in \cite{DN,KoSc,SeW}). Rosenberg proposes a simple
framework extending the classical notion of a Riemannian or Levi-Civita
connection on a Riemannian manifold to the noncommutative tori, proves the
existence and uniqueness of the Levi-Civita connection on the noncommutative
tori, and establishes an analogue of the Gauss-Bonnet theorem for conformal
deformations of the flat metric on the noncommutative tori with a machinery
substantially simpler than that of the highly interesting but technically very
sophisticated spectral geometric approaches to curvature of Connes, Moscovici,
Tretkoff, Fathizadeh, Khalkhali, and others \cite{CM,CT,FKh1,FKh2}.

Intrigued by Rosenberg's paper, we made some observations which seem worthy of
reporting. In particular, in Section 1 of this note, we describe the affine
space of general torsion-free metric-compatible connections (without requiring
that the connection operators for the inner *-derivations act by left
multiplication) on the noncommutative tori and characterize the curvature of
such connections. In Section 2, we find that a Gauss-Bonnet theorem still
holds for a large class of non-commuting (and hence non-conformal) scalings
along the principal directions of a noncommutative two-torus and also for a
class of symplectic scalings.

Our interest in this work was born out of attending a series of conferences
and workshops at the Fields Institute in Toronto, 2013--2015, beginning with
the Focus Program on Noncommutative Geometry and Quantum Groups in celebration
of Marc A. Rieffel's 75th birthday, where Rosenberg introduced his findings on
noncommutative Levi-Civita connections. Both authors would like to thank the
Fields Institute and the organizers of these events, and acknowledge the
support by Jonathan Rosenberg's conference grant NSF-1266158. The second
author also thanks the Mathematics Institute of Academia Sinica, where the
results of this paper were first presented in a colloquium, for the warm
hospitality and support during his visit in the summer of 2015.

\section{Levi-Civita connections}

Let $A_{\Theta}$ be an $n$-dimensional noncommutative torus, that is, the
universal C*-algebra generated by $n$ unitaries $u_{1},...,u_{n}$ satisfying
the relations $u_{j}u_{k}=e^{2\pi i\Theta_{jk}}u_{k}u_{j}$ where $\Theta$ is
an $n\times n$ real skew-symmetric matrix.

We will now recall the definition of Levi-Civita connections for
noncommutative $n$-tori $A_{\Theta}$ given by Rosenberg \cite{Ro}. Throughout
this paper, we work in Rosenberg's framework (except we will remove one
axiom), in particular assuming\ $\Theta$ \textquotedblleft generically
transcendental\textquotedblright\ so that $A_{\Theta}$ is simple and a
classification result of Bratteli, Elliott, and Jorgensen \cite{BEJ} reduces
to the statement that the Lie algebra $\mathcal{D}_{\Theta}$ of all
*-derivations $\delta:A_{\Theta}^{\infty}\rightarrow A_{\Theta}^{\infty}$ is
the $\mathbb{R}$-linear span of the $\partial_{j}$'s and all (bounded) inner
*-derivations, where $A_{\Theta}^{\infty}$ is the subalgebra of $A_{\Theta}$
consisting of all smooth vectors and the $\partial_{j}$'s are the partial
derivatives associated with the canonical $\mathbb{T}^{n}$-action on
$A_{\Theta}$.

Following the general spirit of \cite{C1,C2}, we define a connection for a
noncommutative $n$-torus $A_{\Theta}$ as an $\mathbb{R}$-bilinear map
$\nabla:\mathcal{D}_{\Theta}\times\chi_{\Theta}\rightarrow\chi_{\Theta}$ for
the free left $A_{\Theta}^{\infty}$-module $\chi_{\Theta}:=\left(  A_{\Theta
}^{\infty}\right)  ^{n}$ (which plays the role of the noncommutative tangent
bundle of $A_{\Theta}$ as adopted in \cite{Ro}) such that $\nabla\left(
X,\cdot\right)  \equiv\nabla_{X}$ satisfies
\[
\nabla_{X}\left(  aY\right)  =\left(  X\cdot a\right)  Y+a\nabla_{X}\left(
Y\right)  \text{\ \ for\ \ }\left(  a,Y\right)  \in A_{\Theta}^{\infty}%
\times\chi_{\Theta}.
\]
Fixing a Riemannian metric or more precisely an $A_{\Theta}^{\infty}$-valued
inner product $\left\langle ,\right\rangle $ on $\chi_{\Theta}$ (specified by
a positive invertible matrix $\left(  g_{jk}\right)  \in M_{n}\left(
A_{\Theta}^{\infty}\right)  $ with self-adjoint entries $g_{jk}=\left\langle
\partial_{j},\partial_{k}\right\rangle \in\left(  A_{\Theta}^{\infty}\right)
_{sa}$), a connection $\nabla$ is called a Levi-Civita connection in \cite{Ro}
if it is (1) compatible with the metric, i.e. $X\cdot\left\langle
Y,Z\right\rangle =\left\langle \nabla_{X}Y,Z\right\rangle +\left\langle
Y,\nabla_{X}Z\right\rangle $ for all $X\in\mathcal{D}_{\Theta}$ and
$Y,Z\in\chi_{\Theta}$, (2) torsion-free, i.e. $\nabla_{\partial_{j}}%
\partial_{k}=\nabla_{\partial_{k}}\partial_{j}$ for all $j,k$, (3) real, i.e.
$\left\langle \nabla_{\partial_{j}}\partial_{k},\partial_{l}\right\rangle $ is
self-adjoint in $A_{\Theta}^{\infty}$ for all $j,k,l$, and (4) $\nabla
_{\operatorname{ad}_{\tilde{a}}}=\tilde{a}\cdot$ is left multiplication by
$\tilde{a}$ for any inner *-derivation $\operatorname{ad}_{\tilde{a}}$ given
by a skew-adjoint $\tilde{a}\in A_{\Theta}^{\infty}$ with $\tau\left(
\tilde{a}\right)  =0$ for the unique $\mathbb{T}^{n}$-invariant tracial state
$\tau$ of $A_{\Theta}$. Rosenberg establishes that such a Levi-Civita
connection exists and is unique for any Riemannian metric on $\chi_{\Theta}$.
Among all the notions and conditions mentioned above, only the condition (4)
is not suggested by a corresponding classical notion or condition. So in the
first part of this note, we look into what happens when this condition is
dropped from the properties imposed on Levi-Civita connections.

As in the classical situation, the collection $\mathcal{C}_{\Theta}$ of all
connections $\nabla:\mathcal{D}_{\Theta}\times\chi_{\Theta}\rightarrow
\chi_{\Theta}$ (which are canonically identified as certain $\mathbb{R}%
$-linear maps $\mathcal{D}_{\Theta}\rightarrow\operatorname{End}_{\mathbb{C}%
}\left(  \chi_{\Theta}\right)  $) is an affine space of linear maps
$\mathcal{D}_{\Theta}\rightarrow\operatorname{End}_{A_{\Theta}^{\infty}%
}\left(  \chi_{\Theta}\right)  $, i.e. $\nabla-\nabla^{\prime}\in
\operatorname{Hom}_{\mathbb{R}}\left(  \mathcal{D}_{\Theta},\operatorname{End}%
_{A_{\Theta}^{\infty}}\left(  \chi_{\Theta}\right)  \right)  $ for any
$\nabla,\nabla^{\prime}\in\mathcal{C}_{\Theta}$. Here we use
$\operatorname{Hom}_{\mathbb{R}}$ instead of $\operatorname{Hom}_{A_{\Theta
}^{\infty}}$ since as pointed out in Rosenberg's paper, $\mathcal{D}_{\Theta}$
is no longer a module over the algebra $A_{\Theta}^{\infty}$ of the base
manifold because of the non-commutativity of $A_{\Theta}^{\infty}$.

The inner derivations $\operatorname{ad}_{a}:b\mapsto\left[  a,b\right]
\equiv ab-ba$ of $A_{\Theta}^{\infty}$ with $a\in A_{\Theta}^{\infty}$ are in
one-to-one correspondence with the equivalence classes $\left[  a\right]  $ in
the quotient *-space $A_{\Theta}^{\infty}/\mathbb{C}$, i.e. $\operatorname{ad}%
_{a}=\operatorname{ad}_{a^{\prime}}$ or equivalently $\left[  a-a^{\prime
},b\right]  =0$ for all $b\in A_{\Theta}^{\infty}$ if and only if
$a-a^{\prime}\in\mathbb{C}$, because $\mathbb{C}$ is the center of $A_{\Theta
}$ (and hence of the dense unital subalgebra $A_{\Theta}^{\infty}$) by the
simplicity of the algebra $A_{\Theta}$. Any equivalence class $\left[
a\right]  \in A_{\Theta}^{\infty}/\mathbb{C}$ has a unique representative
$a-\tau\left(  a\right)  \in\ker\left(  \tau\right)  $ for the canonical trace
$\tau$ of $A_{\Theta}$. Furthermore $\operatorname{ad}_{a}$ is a *-derivation
of $A_{\Theta}^{\infty}$ if and only if $\left[  a\right]  $ is skew-adjoint,
i.e. $\left[  a\right]  ^{\ast}\equiv\left[  a^{\ast}\right]  =-\left[
a\right]  $ in $A_{\Theta}^{\infty}/\mathbb{C}$ or $a+a^{\ast}\in\mathbb{C}$,
because $\left[  a,b^{\ast}\right]  =\left[  a,b\right]  ^{\ast}$, or
equivalently $b^{\ast}\left(  a+a^{\ast}\right)  =\left(  a+a^{\ast}\right)
b^{\ast}$ for all $b\in A_{\Theta}^{\infty}$ if and only if $a+a^{\ast}%
\in\mathbb{C}$.

Now observe that the Lie algebra of inner *-derivations of $A_{\Theta}%
^{\infty}$ can be identified with the Lie algebra $i\left(  A_{\Theta}%
^{\infty}\right)  _{sa}\cap\ker\left(  \tau\right)  $ consisting of those
skew-adjoint elements $\tilde{a}\in A_{\Theta}^{\infty}$ with $\tau\left(
\tilde{a}\right)  =0$, where $\left(  A_{\Theta}^{\infty}\right)  _{sa}$
denotes the space of self-adjoint elements of $A_{\Theta}^{\infty}$. Indeed
the inner *-derivations of $A_{\Theta}^{\infty}$ have unique representatives
$a-\tau\left(  a\right)  $ with $a^{\ast}=-a$ modulo $\mathbb{C}$, say,
$a^{\ast}=-a+\left(  x+iy\right)  $ with $x,y\in\mathbb{R}$ such that
$\overline{\tau\left(  a\right)  }=\tau\left(  a^{\ast}\right)  =-\tau\left(
a\right)  +\left(  x+iy\right)  $, or equivalently, $x=2\operatorname{Re}%
\left(  \tau\left(  a\right)  \right)  $ and $y=0$. Now the representative
$\tilde{a}:=a-\tau\left(  a\right)  \in\ker\left(  \tau\right)  $ satisfies%
\[
\tilde{a}^{\ast}=a^{\ast}-\overline{\tau\left(  a\right)  }=\left(
-a+2\operatorname{Re}\left(  \tau\left(  a\right)  \right)  \right)
-\overline{\tau\left(  a\right)  }=-a+\tau\left(  a\right)  =-\tilde{a}%
\]
i.e. $\tilde{a}$ is skew-adjoint. Conversely any skew-adjoint element
$\tilde{a}\in A_{\Theta}^{\infty}$ with $\tau\left(  \tilde{a}\right)  =0$
clearly determines a unique inner *-derivation of $A_{\Theta}^{\infty}$.

Rosenberg's theorem on the existence and uniqueness of the Levi-Civita
connection on the noncommutative tori $A_{\Theta}^{\infty}$ is proved under a
setting (i.e. under condition (4) above) that prescribes and fixes a principal
representative of $\nabla_{\operatorname{ad}_{a}}$ for all inner *-derivations
$\operatorname{ad}_{a}$, namely, $\nabla_{\operatorname{ad}_{a}}\left(
Y\right)  :=\left(  a-\tau\left(  a\right)  \right)  Y$ for $Y\in\chi_{\Theta
}$. It is of interest to ask whether the theorem still holds when more general
connections $\nabla$ for the inner *-derivations $\operatorname{ad}_{a}$ are
allowed (i.e. when only conditions (1)-(3) are required). Since the connection
operators $\nabla_{D}$ for $D$ in the $\mathbb{R}$-span of the $\partial_{i}%
$'s are considered in their most general form in Rosenberg's theorem, we only
need to concentrate on $\nabla_{\operatorname{ad}_{a}}$. More precisely,
letting $\tilde{\nabla}$ denote the unique Levi-Civita connection established
by Rosenberg's theorem for any fixed Riemannian metric $\gamma:=\left(
g_{jk}\right)  \in M_{n}\left(  A_{\Theta}^{\infty}\right)  $ on $\chi
_{\Theta}$ with self-adjoint $g_{jk}=\left\langle \partial_{j},\partial
_{k}\right\rangle \in\left(  A_{\Theta}^{\infty}\right)  _{sa}$, we consider
all connections $\nabla$ of the form $\nabla_{D}=\tilde{\nabla}_{D}$ for $D$
in the $\mathbb{R}$-span of $\partial_{i}$ and $\nabla_{\operatorname{ad}%
_{\tilde{a}}}=\tilde{\nabla}_{\operatorname{ad}_{\tilde{a}}}+\mu\left(
\tilde{a}\right)  $ for some
\[
\mu\in\operatorname{Hom}_{\mathbb{R}}\left(  i\left(  A_{\Theta}^{\infty
}\right)  _{sa}\cap\ker\left(  \tau\right)  ,\operatorname{End}_{A_{\Theta
}^{\infty}}\left(  \chi_{\Theta}\right)  \right)  ,
\]
i.e. $\nabla_{\operatorname{ad}_{\tilde{a}}}\left(  Y\right)  =\tilde{a}%
Y+\mu\left(  \tilde{a}\right)  \left(  Y\right)  $ for $Y\in\chi_{\Theta}$,
where $\tilde{a}\in i\left(  A_{\Theta}^{\infty}\right)  _{sa}\cap\ker\left(
\tau\right)  $. Note that since the torsion-free constraint involves only the
$\partial_{j}$'s, such a connection is automatically torsion-free. We ask what
condition on $\mu$ makes such a connection $\nabla$ compatible with the fixed
Riemannian metric $\gamma$ on $\chi_{\Theta}$, and whether it makes $\nabla$ unique.

We canonically identify an endomorphism $\mu\left(  \tilde{a}\right)  $ of the
free $A_{\Theta}^{\infty}$-module $\chi_{\Theta}$ with an element $\mu
_{\tilde{a}}$ of $M_{n}\left(  A_{\Theta}^{\infty}\right)  $. Then since
\[
\operatorname{ad}_{\tilde{a}}\left(  \left\langle X,Y\right\rangle \right)
=\left\langle \tilde{\nabla}_{\operatorname{ad}_{\tilde{a}}}X,Y\right\rangle
+\left\langle X,\tilde{\nabla}_{\operatorname{ad}_{\tilde{a}}}Y\right\rangle
\]
holds for the Levi-Civita connection $\tilde{\nabla}$, the
metric-compatibility condition
\[
\operatorname{ad}_{\tilde{a}}\left(  \left\langle X,Y\right\rangle \right)
=\left\langle \nabla_{\operatorname{ad}_{\tilde{a}}}X,Y\right\rangle
+\left\langle X,\nabla_{\operatorname{ad}_{\tilde{a}}}Y\right\rangle
\]
is equivalent to%
\[
0=\left\langle \mu\left(  \tilde{a}\right)  \left(  X\right)  ,Y\right\rangle
+\left\langle X,\mu\left(  \tilde{a}\right)  \left(  Y\right)  \right\rangle
\]
i.e. $0=\mu_{\tilde{a}}\gamma+\gamma\left(  \mu_{\tilde{a}}\right)  ^{\ast}$,
or equivalently,
\[
\left(  \mu_{\tilde{a}}\right)  ^{\ast}=-\gamma^{-1}\mu_{\tilde{a}}%
\gamma\equiv-\operatorname{Ad}_{\gamma}^{-1}\left(  \mu_{\tilde{a}}\right)  ,
\]
or $\mu_{\tilde{a}}\in\ker\left(  ^{\ast}+\operatorname{Ad}_{\gamma}%
^{-1}\right)  \subset M_{n}\left(  A_{\Theta}^{\infty}\right)  $.

Any skew-adjoint element $\nu\in M_{n}\left(  A_{\Theta}^{\infty}\right)  $
that commutes with $\gamma$ satisfies $\nu^{\ast}=-\gamma^{-1}\nu\gamma$, for
example, any $\nu=if\left(  \gamma\right)  $ for a real polynomial (or
analytic) $f$ on $\sigma\left(  \gamma\right)  \subset\left(  0,\infty\right)
$. So any nonzero $\mathbb{R}$-linear map $\mu$ from $i\left(  A_{\Theta
}^{\infty}\right)  _{sa}\cap\ker\left(  \tau\right)  $ into the space of all
skew-adjoint elements $\nu\in M_{n}\left(  A_{\Theta}^{\infty}\right)  $ in
the commutant of $\gamma$ gives rise to a torsion-free connection $\nabla$
compatible with the fixed Riemannian metric $\gamma$. Such a map $\mu$ (and
hence $\nabla$) is highly non-unique. We summarize this observation as follows:

\textbf{Proposition 1}. The space of all torsion-free connections on
$\chi_{\Theta}$ compatible with a Riemannian metric given by an invertible
positive matrix $\gamma\in M_{n}\left(  A_{\Theta}^{\infty}\right)  $ is the
nontrivial affine space
\[
\tilde{\nabla}+\operatorname{Hom}_{\mathbb{R}}\left(  \mathfrak{g},\ker\left(
^{\ast}+\operatorname{Ad}_{\gamma}^{-1}\right)  \right)  ,
\]
where $\mathfrak{g}:=i\left(  A_{\Theta}^{\infty}\right)  _{sa}\cap\ker\left(
\tau\right)  $ is identified with the Lie algebra of inner *-derivations of
$A_{\Theta}^{\infty}$ under the adjoint representation $\operatorname{ad}$,
$\tilde{\nabla}$ is Rosenberg's noncommutative Levi-Civita connection
determined by $\gamma$, and $\tilde{\nabla}+\mu$ is set to equal
$\tilde{\nabla}$ on the $\mathbb{R}$-span of the $\partial_{j}$'s for any
$\mu\in\operatorname{Hom}_{\mathbb{R}}\left(  \mathfrak{g},\ker\left(  ^{\ast
}+\operatorname{Ad}_{\gamma}^{-1}\right)  \right)  $.%

{\bigskip}%

When only the connections with $\nabla_{\operatorname{ad}_{\tilde{a}}}\left(
Z\right)  =\tilde{a}Z$ for $Z\in\chi_{\Theta}$ are considered, Rosenberg shows
that $\nabla_{\operatorname{ad}_{\tilde{a}}}$ does not play a role in the
curvature form \cite{Ro}
\[
R\left(  X,Y\right)  :=\nabla_{Y}\nabla_{X}-\nabla_{X}\nabla_{Y}%
+\nabla_{\left[  X,Y\right]  }\in\operatorname{End}_{A_{\Theta}^{\infty}%
}\left(  \chi_{\Theta}\right)
\]
for $X,Y\in\mathcal{D}_{\Theta}$, i.e. $R\left(  \operatorname{ad}_{\tilde{a}%
},X\right)  =0$ for all $X\in\mathcal{D}_{\Theta}$.

However for more general connections $\nabla$ with $\nabla_{\operatorname{ad}%
_{\tilde{a}}}=\tilde{a}\cdot+\mu\left(  \tilde{a}\right)  $, we can get
non-vanishing $R\left(  \operatorname{ad}_{\tilde{a}},X\right)  $. Indeed with
$\left[  \operatorname{ad}_{\tilde{a}},X\right]  =\operatorname{ad}%
_{-X\cdot\tilde{a}}$ in $\mathcal{D}_{\Theta}$, for any $b\in A_{\Theta
}^{\infty}$,
\[
R\left(  \operatorname{ad}_{\tilde{a}},X\right)  \left(  b\partial_{k}\right)
=\nabla_{X}\left(  \left(  \tilde{a}\cdot+\mu\left(  \tilde{a}\right)
\right)  \left(  b\partial_{k}\right)  \right)  -\left(  \tilde{a}\cdot
+\mu\left(  \tilde{a}\right)  \right)  \left(  \nabla_{X}\left(  b\partial
_{k}\right)  \right)  +\nabla_{\operatorname{ad}_{-X\cdot\tilde{a}}}\left(
b\partial_{k}\right)
\]%
\[
=\nabla_{X}\left(  \tilde{a}b\partial_{k}+b\mu\left(  \tilde{a}\right)
\left(  \partial_{k}\right)  \right)  -\left(  \tilde{a}\cdot+\mu\left(
\tilde{a}\right)  \right)  \left(  \left(  X\cdot b\right)  \partial
_{k}+b\nabla_{X}\left(  \partial_{k}\right)  \right)  -\left(  X\cdot\tilde
{a}\right)  b\partial_{k}-b\mu\left(  X\cdot\tilde{a}\right)  \left(
\partial_{k}\right)
\]%
\begin{align*}
&  =\left(  X\cdot\left(  \tilde{a}b\right)  \right)  \partial_{k}+\tilde
{a}b\nabla_{X}\left(  \partial_{k}\right)  +\left(  X\cdot b\right)
\mu\left(  \tilde{a}\right)  \left(  \partial_{k}\right)  +b\nabla_{X}\left(
\mu\left(  \tilde{a}\right)  \left(  \partial_{k}\right)  \right)  -\tilde
{a}\left(  X\cdot b\right)  \partial_{k}\\
&  -\left(  X\cdot b\right)  \mu\left(  \tilde{a}\right)  \left(  \partial
_{k}\right)  -\tilde{a}b\nabla_{X}\left(  \partial_{k}\right)  -b\mu\left(
\tilde{a}\right)  \left(  \nabla_{X}\left(  \partial_{k}\right)  \right)
-\left(  X\cdot\tilde{a}\right)  b\partial_{k}-b\mu\left(  X\cdot\tilde
{a}\right)  \left(  \partial_{k}\right)
\end{align*}%
\[
=b\nabla_{X}\left(  \mu\left(  \tilde{a}\right)  \left(  \partial_{k}\right)
\right)  -b\mu\left(  \tilde{a}\right)  \left(  \nabla_{X}\left(  \partial
_{k}\right)  \right)  -b\mu\left(  X\cdot\tilde{a}\right)  \left(
\partial_{k}\right)
\]%
\[
=b\left(  \left[  \nabla_{X},\mu\left(  \tilde{a}\right)  \right]  -\mu\left(
X\cdot\tilde{a}\right)  \right)  \left(  \partial_{k}\right)
\]
which shows that $R\left(  \operatorname{ad}_{\tilde{a}},X\right)  =\left[
\nabla_{X},\mu\left(  \tilde{a}\right)  \right]  -\mu\left(  X\cdot\tilde
{a}\right)  $, where $X\cdot\tilde{a}\in\mathfrak{g}$ for $\tilde{a}%
\in\mathfrak{g}$, and $\left[  \nabla_{X},\mu\left(  \tilde{a}\right)
\right]  \in\operatorname{End}_{A_{\Theta}^{\infty}}\left(  \chi_{\Theta
}\right)  $ is the commutator bracket of $\nabla_{X}$ and $\mu\left(
\tilde{a}\right)  $ in $\operatorname{End}_{\mathbb{C}}\left(  \chi_{\Theta
}\right)  $. So $R\left(  \operatorname{ad}_{\tilde{a}},X\right)  =0$ if and
only if $\left[  \nabla_{X},\mu\left(  \tilde{a}\right)  \right]  =\mu\left(
X\cdot\tilde{a}\right)  $, which is a nontrivial condition on $\mu$. Actually
for $X=\operatorname{ad}_{\tilde{b}}$, this condition can be rewritten as
$\left[  \mu\left(  \tilde{b}\right)  ,\mu\left(  \tilde{a}\right)  \right]
=\mu\left(  \left[  \tilde{b},\tilde{a}\right]  \right)  $. So $R\left(
\operatorname{ad}_{\tilde{a}},\operatorname{ad}_{\tilde{b}}\right)  =0$ for
all inner *-derivations $\operatorname{ad}_{\tilde{a}},\operatorname{ad}%
_{\tilde{b}}$ if and only if $\mu$ is a Lie algebra homomorphism. We summarize
these observations as follows:

\textbf{Proposition 2}. For a connection $\nabla$ on $\chi_{\Theta}$ with
$\nabla_{\operatorname{ad}_{\tilde{a}}}=\tilde{a}\cdot+\mu\left(  \tilde
{a}\right)  $ for some $\mu\in\operatorname{Hom}_{\mathbb{R}}\left(
\mathfrak{g},\operatorname{End}_{A_{\Theta}^{\infty}}\left(  \chi_{\Theta
}\right)  \right)  $, the curvature form $R\left(  \operatorname{ad}%
_{\tilde{a}},X\right)  $ equals
\[
\left[  \nabla_{X},\mu\right]  \left(  \tilde{a}\right)  :=\left[  \nabla
_{X},\mu\left(  \tilde{a}\right)  \right]  -\mu\left(  X\cdot\tilde{a}\right)
\]
for $X\in\mathcal{D}_{\Theta}$ and inner *-derivations $\tilde{a}%
\in\mathfrak{g}:=i\left(  A_{\Theta}^{\infty}\right)  _{sa}\cap\ker\left(
\tau\right)  $, and $R$ vanishes completely on all pairs of inner
*-derivations if and only if $\mu$ is a Lie algebra homomorphism.

\section{Gauss-Bonnet theorem}

A major result of Rosenberg \cite{Ro} is that the Gauss-Bonnet theorem holds
for the noncommutative two-tori $A_{\theta}$ when $\chi_{\theta}\equiv\left(
A_{\theta}^{\infty}\right)  ^{2}$ is equipped with a conformal deformation of
the flat metric, i.e. $\left(  g_{jk}\right)  =e^{h}I_{2}$ for $h\in
A_{\theta}^{\infty}$ self-adjoint, in the sense that the Riemannian integral
$\tau\left(  R_{1,2,1,2}e^{-h}\right)  $ of the Gaussian curvature of the
associated Levi-Civita connection vanishes as the Euler number of
$\mathbb{T}^{2}$, where $R_{1,2,1,2}:=\left\langle R\left(  \partial
_{1},\partial_{2}\right)  \partial_{1},\partial_{2}\right\rangle $, and
$e^{-h}=e^{-2h}e^{h}$ reflects the normalization of the non-unit vectors
$\partial_{1},\partial_{2}$ and the deformed volume element of $\mathbb{T}%
^{2}$. In the following, we show that for two broad classes of non-conformal
deformations of the flat metric, such a Gauss-Bonnet theorem still holds.

Indeed we consider the class of \textquotedblleft diagonal\textquotedblright%
\ non-conformal metrics of the form $\left\langle \partial_{j},\partial
_{k}\right\rangle \equiv g_{jk}=\delta_{jk}a_{k}$ for some distinct positive
invertible elements $a_{1},a_{2}\in A_{\theta}^{\infty}$, and show that
$\tau\left(  a_{1}^{\frac{-1}{2}}R_{1,2,1,2}a_{2}^{\frac{-1}{2}}\right)  $
vanishes either when each $a_{j}$ commutes with its first-order derivatives
$\partial_{1}a_{j}$ and $\partial_{2}a_{j}$ or when $a_{1}a_{2}=1$.

We first compute the Levi-Civita connection for the \textquotedblleft
non-conformal\textquotedblright\ metric
\[
\gamma=\left(  g_{jk}\right)  =\left(
\begin{array}
[c]{cc}%
a_{1} & 0\\
0 & a_{2}%
\end{array}
\right)
\]
for arbitrary positive invertible elements $a_{1},a_{2}\in A_{\theta}^{\infty
}$. From the formula
\[
\left\langle \nabla_{j}\partial_{k},\partial_{l}\right\rangle =\frac{1}%
{2}\left(  \partial_{j}\left\langle \partial_{k},\partial_{l}\right\rangle
+\partial_{k}\left\langle \partial_{j},\partial_{l}\right\rangle -\partial
_{l}\left\langle \partial_{j},\partial_{k}\right\rangle \right)  =\frac{1}%
{2}\left(  \delta_{kl}\partial_{j}a_{k}+\delta_{jl}\partial_{k}a_{j}%
-\delta_{jk}\partial_{l}a_{j}\right)
\]
established by Rosenberg where $\nabla_{j}:=\nabla_{\partial_{j}}$, we get%
\[
\left(  \left\langle \nabla_{1}\partial_{k},\partial_{l}\right\rangle \right)
_{kl}=\frac{1}{2}\left(
\begin{array}
[c]{cc}%
\partial_{1}a_{1} & -\partial_{2}a_{1}\\
\partial_{2}a_{1} & \partial_{1}a_{2}%
\end{array}
\right)
\]
and%
\[
\left(  \left\langle \nabla_{2}\partial_{k},\partial_{l}\right\rangle \right)
_{kl}=\frac{1}{2}\left(
\begin{array}
[c]{cc}%
\partial_{2}a_{1} & \partial_{1}a_{2}\\
-\partial_{1}a_{2} & \partial_{2}a_{2}%
\end{array}
\right)  ,
\]
which imply the identities%
\[
\left\{
\begin{array}
[c]{l}%
\nabla_{1}\partial_{1}=\frac{1}{2}\left(  \left(  \partial_{1}a_{1}\right)
a_{1}^{-1}\partial_{1}-\left(  \partial_{2}a_{1}\right)  a_{2}^{-1}%
\partial_{2}\right)  ,\\
\nabla_{1}\partial_{2}=\frac{1}{2}\left(  \left(  \partial_{2}a_{1}\right)
a_{1}^{-1}\partial_{1}+\left(  \partial_{1}a_{2}\right)  a_{2}^{-1}%
\partial_{2}\right)  ,\\
\nabla_{2}\partial_{1}=\frac{1}{2}\left(  \left(  \partial_{2}a_{1}\right)
a_{1}^{-1}\partial_{1}+\left(  \partial_{1}a_{2}\right)  a_{2}^{-1}%
\partial_{2}\right)  ,\\
\nabla_{2}\partial_{2}=\frac{1}{2}\left(  -\left(  \partial_{1}a_{2}\right)
a_{1}^{-1}\partial_{1}+\left(  \partial_{2}a_{2}\right)  a_{2}^{-1}%
\partial_{2}\right)  .
\end{array}
\right.
\]
Thus
\[
4R\left(  \partial_{1},\partial_{2}\right)  \left(  \partial_{1}\right)
=4\left(  \nabla_{2}\nabla_{1}-\nabla_{1}\nabla_{2}\right)  \left(
\partial_{1}\right)
\]%
\[
=2\nabla_{2}\left(  \left(  \partial_{1}a_{1}\right)  a_{1}^{-1}\partial
_{1}-\left(  \partial_{2}a_{1}\right)  a_{2}^{-1}\partial_{2}\right)
-2\nabla_{1}\left(  \left(  \partial_{2}a_{1}\right)  a_{1}^{-1}\partial
_{1}+\left(  \partial_{1}a_{2}\right)  a_{2}^{-1}\partial_{2}\right)
\]%
\begin{align*}
&  =2\partial_{2}\left(  \left(  \partial_{1}a_{1}\right)  a_{1}^{-1}\right)
\partial_{1}+\left(  \partial_{1}a_{1}\right)  a_{1}^{-1}\left(  \left(
\partial_{2}a_{1}\right)  a_{1}^{-1}\partial_{1}+\left(  \partial_{1}%
a_{2}\right)  a_{2}^{-1}\partial_{2}\right) \\
&  -2\partial_{2}\left(  \left(  \partial_{2}a_{1}\right)  a_{2}^{-1}\right)
\partial_{2}-\left(  \partial_{2}a_{1}\right)  a_{2}^{-1}\left(  -\left(
\partial_{1}a_{2}\right)  a_{1}^{-1}\partial_{1}+\left(  \partial_{2}%
a_{2}\right)  a_{2}^{-1}\partial_{2}\right) \\
&  -2\partial_{1}\left(  \left(  \partial_{2}a_{1}\right)  a_{1}^{-1}\right)
\partial_{1}-\left(  \partial_{2}a_{1}\right)  a_{1}^{-1}\left(  \left(
\partial_{1}a_{1}\right)  a_{1}^{-1}\partial_{1}-\left(  \partial_{2}%
a_{1}\right)  a_{2}^{-1}\partial_{2}\right) \\
&  -2\partial_{1}\left(  \left(  \partial_{1}a_{2}\right)  a_{2}^{-1}\right)
\partial_{2}-\left(  \partial_{1}a_{2}\right)  a_{2}^{-1}\left(  \left(
\partial_{2}a_{1}\right)  a_{1}^{-1}\partial_{1}+\left(  \partial_{1}%
a_{2}\right)  a_{2}^{-1}\partial_{2}\right)
\end{align*}%
\begin{align*}
&  =2\partial_{2}\left(  \left(  \partial_{1}a_{1}\right)  a_{1}^{-1}\right)
\partial_{1}+\left(  \partial_{1}a_{1}\right)  a_{1}^{-1}\left(  \partial
_{2}a_{1}\right)  a_{1}^{-1}\partial_{1}+\left(  \partial_{1}a_{1}\right)
a_{1}^{-1}\left(  \partial_{1}a_{2}\right)  a_{2}^{-1}\partial_{2}\\
&  -2\partial_{2}\left(  \left(  \partial_{2}a_{1}\right)  a_{2}^{-1}\right)
\partial_{2}+\left(  \partial_{2}a_{1}\right)  a_{2}^{-1}\left(  \partial
_{1}a_{2}\right)  a_{1}^{-1}\partial_{1}-\left(  \partial_{2}a_{1}\right)
a_{2}^{-1}\left(  \partial_{2}a_{2}\right)  a_{2}^{-1}\partial_{2}\\
&  -2\partial_{1}\left(  \left(  \partial_{2}a_{1}\right)  a_{1}^{-1}\right)
\partial_{1}-\left(  \partial_{2}a_{1}\right)  a_{1}^{-1}\left(  \partial
_{1}a_{1}\right)  a_{1}^{-1}\partial_{1}+\left(  \partial_{2}a_{1}\right)
a_{1}^{-1}\left(  \partial_{2}a_{1}\right)  a_{2}^{-1}\partial_{2}\\
&  -2\partial_{1}\left(  \left(  \partial_{1}a_{2}\right)  a_{2}^{-1}\right)
\partial_{2}-\left(  \partial_{1}a_{2}\right)  a_{2}^{-1}\left(  \partial
_{2}a_{1}\right)  a_{1}^{-1}\partial_{1}-\left(  \partial_{1}a_{2}\right)
a_{2}^{-1}\left(  \partial_{1}a_{2}\right)  a_{2}^{-1}\partial_{2}%
\end{align*}%
\begin{align*}
&  \overset{\text{**}}{=}2\left(  \partial_{2}\partial_{1}a_{1}\right)
a_{1}^{-1}\partial_{1}-\left(  \partial_{1}a_{1}\right)  a_{1}^{-1}\left(
\partial_{2}a_{1}\right)  a_{1}^{-1}\partial_{1}+\left(  \partial_{1}%
a_{1}\right)  a_{1}^{-1}\left(  \partial_{1}a_{2}\right)  a_{2}^{-1}%
\partial_{2}\\
&  -2\left(  \partial_{2}\partial_{2}a_{1}\right)  a_{2}^{-1}\partial
_{2}+\left(  \partial_{2}a_{1}\right)  a_{2}^{-1}\left(  \partial_{1}%
a_{2}\right)  a_{1}^{-1}\partial_{1}+\left(  \partial_{2}a_{1}\right)
a_{2}^{-1}\left(  \partial_{2}a_{2}\right)  a_{2}^{-1}\partial_{2}\\
&  -2\left(  \partial_{1}\partial_{2}a_{1}\right)  a_{1}^{-1}\partial
_{1}+\left(  \partial_{2}a_{1}\right)  a_{1}^{-1}\left(  \partial_{1}%
a_{1}\right)  a_{1}^{-1}\partial_{1}+\left(  \partial_{2}a_{1}\right)
a_{1}^{-1}\left(  \partial_{2}a_{1}\right)  a_{2}^{-1}\partial_{2}\\
&  -2\left(  \partial_{1}\partial_{1}a_{2}\right)  a_{2}^{-1}\partial
_{2}-\left(  \partial_{1}a_{2}\right)  a_{2}^{-1}\left(  \partial_{2}%
a_{1}\right)  a_{1}^{-1}\partial_{1}+\left(  \partial_{1}a_{2}\right)
a_{2}^{-1}\left(  \partial_{1}a_{2}\right)  a_{2}^{-1}\partial_{2}%
\end{align*}%
\begin{align*}
&  =\left(  -\left(  \partial_{1}a_{1}\right)  a_{1}^{-1}\left(  \partial
_{2}a_{1}\right)  +\left(  \partial_{2}a_{1}\right)  a_{2}^{-1}\left(
\partial_{1}a_{2}\right)  +\left(  \partial_{2}a_{1}\right)  a_{1}^{-1}\left(
\partial_{1}a_{1}\right)  -\left(  \partial_{1}a_{2}\right)  a_{2}^{-1}\left(
\partial_{2}a_{1}\right)  \right)  a_{1}^{-1}\partial_{1}\\
&  +\left(  \left(  \partial_{1}a_{1}\right)  a_{1}^{-1}\left(  \partial
_{1}a_{2}\right)  +\left(  \partial_{2}a_{1}\right)  a_{2}^{-1}\left(
\partial_{2}a_{2}\right)  +\left(  \partial_{2}a_{1}\right)  a_{1}^{-1}\left(
\partial_{2}a_{1}\right)  +\left(  \partial_{1}a_{2}\right)  a_{2}^{-1}\left(
\partial_{1}a_{2}\right)  \right)  a_{2}^{-1}\partial_{2}\\
&  +2\left(  \partial_{2}\partial_{1}a_{1}-\partial_{1}\partial_{2}%
a_{1}\right)  a_{1}^{-1}\partial_{1}-2\left(  \partial_{1}\partial_{1}%
a_{2}+\partial_{2}\partial_{2}a_{1}\right)  a_{2}^{-1}\partial_{2}\ ,
\end{align*}
where ** is due to $\partial_{j}\left(  a^{-1}\right)  =-a^{-1}\left(
\partial_{j}a\right)  a^{-1}$ for invertible $a$ (since $a\partial_{j}\left(
a^{-1}\right)  +\left(  \partial_{j}a\right)  a^{-1}=\partial_{j}\left(
aa^{-1}\right)  =0$).

Now we have
\[
4R_{1,2,1,2}=\left\langle 4R\left(  \partial_{1},\partial_{2}\right)  \left(
\partial_{1}\right)  ,\partial_{2}\right\rangle =\left\langle 4\left(
\nabla_{2}\nabla_{1}-\nabla_{1}\nabla_{2}\right)  \left(  \partial_{1}\right)
,\partial_{2}\right\rangle
\]%
\[
=\left(  \partial_{1}a_{1}\right)  a_{1}^{-1}\left(  \partial_{1}a_{2}\right)
+\left(  \partial_{2}a_{1}\right)  a_{2}^{-1}\left(  \partial_{2}a_{2}\right)
+\left(  \partial_{2}a_{1}\right)  a_{1}^{-1}\left(  \partial_{2}a_{1}\right)
+\left(  \partial_{1}a_{2}\right)  a_{2}^{-1}\left(  \partial_{1}a_{2}\right)
-2\left(  \partial_{1}\partial_{1}a_{2}+\partial_{2}\partial_{2}a_{1}\right)
.
\]
Since $\partial_{1}$ and $\partial_{2}$ are orthogonal and of
\textquotedblleft noncommutative\textquotedblright\ lengths $\sqrt{a_{1}}$ and
$\sqrt{a_{2}}$, respectively, we need to scale $R_{1,2,1,2}$ to get a
noncommutative Gaussian curvature $K:=a_{1}^{-1}R_{1,2,1,2}a_{2}^{-1}$. Then
we integrate $K$ against the \textquotedblleft Riemannian\textquotedblright%
\ volume form by taking
\[
\tau\left(  \sqrt{a_{1}}K\sqrt{a_{2}}\right)  \equiv\tau\left(  a_{1}%
^{\frac{-1}{2}}R_{1,2,1,2}a_{2}^{\frac{-1}{2}}\right)  .
\]

We now show that when each $a_{j}$ commutes with both $\partial_{1}a_{j}$ and
$\partial_{2}a_{j}$ (but $a_{1}$ may not commute with $a_{2}$, for example, if
$a_{j}\in C^{\ast}\left(  u_{j}\right)  $ are positive invertible elements of
$A_{\theta}^{\infty}$, for the two generating unitaries $u_{1},u_{2}$ of
$A_{\theta}$), the Gauss-Bonnet theorem holds, i.e.%
\[
\tau\left(  \sqrt{a_{1}}K\sqrt{a_{2}}\right)  =0.
\]

Indeed, it is well known that for an analytic function $f$ on the right half
plane $\operatorname{Re}\left(  z\right)  >0$, such as $z\mapsto z^{\frac
{-1}{2}}$, the ordinary differential calculus applies to $f\left(  a\right)
$, for example $\partial_{j}\left(  a^{\frac{-1}{2}}\right)  =\frac{-1}%
{2}a^{\frac{-3}{2}}\left(  \partial_{j}a\right)  $, when a positive invertible
element $a$ in $A_{\theta}^{\infty}$ commutes with its partial derivatives
$\partial_{j}a$. So we can compute
\[
4a_{1}^{\frac{-1}{2}}R_{1,2,1,2}a_{2}^{\frac{-1}{2}}=a_{1}^{\frac{-1}{2}%
}\left(  \partial_{1}a_{1}\right)  a_{1}^{-1}\left(  \partial_{1}a_{2}\right)
a_{2}^{\frac{-1}{2}}+a_{1}^{\frac{-1}{2}}\left(  \partial_{2}a_{1}\right)
a_{2}^{-1}\left(  \partial_{2}a_{2}\right)  a_{2}^{\frac{-1}{2}}+a_{1}%
^{\frac{-1}{2}}\left(  \partial_{2}a_{1}\right)  a_{1}^{-1}\left(
\partial_{2}a_{1}\right)  a_{2}^{\frac{-1}{2}}%
\]%
\[
+a_{1}^{\frac{-1}{2}}\left(  \partial_{1}a_{2}\right)  a_{2}^{-1}\left(
\partial_{1}a_{2}\right)  a_{2}^{\frac{-1}{2}}-2a_{1}^{\frac{-1}{2}}\left(
\partial_{1}\partial_{1}a_{2}\right)  a_{2}^{\frac{-1}{2}}-2a_{1}^{\frac
{-1}{2}}\left(  \partial_{2}\partial_{2}a_{1}\right)  a_{2}^{\frac{-1}{2}}%
\]%
\[
=a_{1}^{\frac{-3}{2}}\left(  \partial_{1}a_{1}\right)  \left(  \partial
_{1}a_{2}\right)  a_{2}^{\frac{-1}{2}}+a_{1}^{\frac{-1}{2}}\left(
\partial_{2}a_{1}\right)  \left(  \partial_{2}a_{2}\right)  a_{2}^{\frac
{-3}{2}}+a_{1}^{\frac{-3}{2}}\left(  \partial_{2}a_{1}\right)  \left(
\partial_{2}a_{1}\right)  a_{2}^{\frac{-1}{2}}%
\]%
\[
+a_{1}^{\frac{-1}{2}}\left(  \partial_{1}a_{2}\right)  \left(  \partial
_{1}a_{2}\right)  a_{2}^{\frac{-3}{2}}-2a_{1}^{\frac{-1}{2}}\left(
\partial_{1}\partial_{1}a_{2}\right)  a_{2}^{\frac{-1}{2}}-2a_{1}^{\frac
{-1}{2}}\left(  \partial_{2}\partial_{2}a_{1}\right)  a_{2}^{\frac{-1}{2}}%
\]%
\[
=-2\partial_{1}\left(  a_{1}^{\frac{-1}{2}}\right)  \left(  \partial_{1}%
a_{2}\right)  a_{2}^{\frac{-1}{2}}-2a_{1}^{\frac{-1}{2}}\left(  \partial
_{2}a_{1}\right)  \partial_{2}\left(  a_{2}^{\frac{-1}{2}}\right)
-2\partial_{2}\left(  a_{1}^{\frac{-1}{2}}\right)  \left(  \partial_{2}%
a_{1}\right)  a_{2}^{\frac{-1}{2}}%
\]%
\[
-2a_{1}^{\frac{-1}{2}}\left(  \partial_{1}a_{2}\right)  \partial_{1}\left(
a_{2}^{\frac{-1}{2}}\right)  -2a_{1}^{\frac{-1}{2}}\left(  \partial
_{1}\partial_{1}a_{2}\right)  a_{2}^{\frac{-1}{2}}-2a_{1}^{\frac{-1}{2}%
}\left(  \partial_{2}\partial_{2}a_{1}\right)  a_{2}^{\frac{-1}{2}}%
\]%
\[
=-2\partial_{1}\left(  a_{1}^{\frac{-1}{2}}\right)  \left(  \partial_{1}%
a_{2}\right)  a_{2}^{\frac{-1}{2}}-2a_{1}^{\frac{-1}{2}}\left(  \partial
_{1}\partial_{1}a_{2}\right)  a_{2}^{\frac{-1}{2}}-2a_{1}^{\frac{-1}{2}%
}\left(  \partial_{1}a_{2}\right)  \partial_{1}\left(  a_{2}^{\frac{-1}{2}%
}\right)
\]%
\[
-2\partial_{2}\left(  a_{1}^{\frac{-1}{2}}\right)  \left(  \partial_{2}%
a_{1}\right)  a_{2}^{\frac{-1}{2}}-2a_{1}^{\frac{-1}{2}}\left(  \partial
_{2}\partial_{2}a_{1}\right)  a_{2}^{\frac{-1}{2}}-2a_{1}^{\frac{-1}{2}%
}\left(  \partial_{2}a_{1}\right)  \partial_{2}\left(  a_{2}^{\frac{-1}{2}%
}\right)
\]%
\[
=-2\partial_{1}\left(  a_{1}^{\frac{-1}{2}}\left(  \partial_{1}a_{2}\right)
a_{2}^{\frac{-1}{2}}\right)  -2\partial_{2}\left(  a_{1}^{\frac{-1}{2}}\left(
\partial_{2}a_{1}\right)  a_{2}^{\frac{-1}{2}}\right)  ,
\]
which shows that
\[
\tau\left(  a_{1}^{\frac{-1}{2}}R_{1,2,1,2}a_{2}^{\frac{-1}{2}}\right)
=\frac{-2}{4}\tau\left(  \partial_{1}\left(  a_{1}^{\frac{-1}{2}}\left(
\partial_{1}a_{2}\right)  a_{2}^{\frac{-1}{2}}\right)  +\partial_{2}\left(
a_{1}^{\frac{-1}{2}}\left(  \partial_{2}a_{1}\right)  a_{2}^{\frac{-1}{2}%
}\right)  \right)  =0.
\]

We also note that when the diagonal metric is a (non-conformal, but
symplectic) deformation of the flat metric with $a_{1}a_{2}=1$, the
Gauss-Bonnet theorem also holds. Indeed
\[
\tau\left(  a_{1}^{\frac{-1}{2}}R_{1,2,1,2}a_{2}^{\frac{-1}{2}}\right)
=\tau\left(  R_{1,2,1,2}a_{2}^{\frac{-1}{2}}a_{1}^{\frac{-1}{2}}\right)
=\tau\left(  R_{1,2,1,2}\right)
\]%
\[
=\frac{\tau}{4}\left(
\begin{array}
[c]{c}%
\partial_{1}\left(  a_{2}^{-1}\right)  a_{2}\left(  \partial_{1}a_{2}\right)
+\partial_{2}\left(  a_{2}^{-1}\right)  a_{2}^{-1}\left(  \partial_{2}%
a_{2}\right)  +\partial_{2}\left(  a_{2}^{-1}\right)  a_{2}\partial_{2}\left(
a_{2}^{-1}\right) \\
+\left(  \partial_{1}a_{2}\right)  a_{2}^{-1}\left(  \partial_{1}a_{2}\right)
-2\left(  \partial_{1}\partial_{1}a_{2}+\partial_{2}\partial_{2}a_{1}\right)
\end{array}
\right)
\]%
\[
=\frac{\tau}{4}\left(  -a_{2}^{-1}\left(  \partial_{1}a_{2}\right)  \left(
\partial_{1}a_{2}\right)  -\partial_{2}\left(  a_{2}^{-1}\right)  \partial
_{2}\left(  a_{2}^{-1}\right)  a_{2}+\partial_{2}\left(  a_{2}^{-1}\right)
a_{2}\partial_{2}\left(  a_{2}^{-1}\right)  +\left(  \partial_{1}a_{2}\right)
a_{2}^{-1}\left(  \partial_{1}a_{2}\right)  \right)
\]%
\[
=\frac{\tau}{4}\left(  -a_{2}^{-1}\left(  \partial_{1}a_{2}\right)  \left(
\partial_{1}a_{2}\right)  -\partial_{2}\left(  a_{2}^{-1}\right)  \partial
_{2}\left(  a_{2}^{-1}\right)  a_{2}+\partial_{2}\left(  a_{2}^{-1}\right)
\partial_{2}\left(  a_{2}^{-1}\right)  a_{2}+a_{2}^{-1}\left(  \partial
_{1}a_{2}\right)  \left(  \partial_{1}a_{2}\right)  \right)  =0.
\]

We summarize these results as follows:

\textbf{Theorem}. For the noncommutative two-torus $A_{\theta}$, if
$\chi_{\theta}\equiv\left(  A_{\theta}^{\infty}\right)  ^{2}$ is equipped with
a diagonal\ metric of the form $\left\langle \partial_{j},\partial
_{k}\right\rangle \equiv g_{jk}=\delta_{jk}a_{k}$ for some positive invertible
elements $a_{1},a_{2}\in A_{\theta}^{\infty}$, then the Riemannian integral of
the associated Gaussian curvature vanishes, i.e.
\[
\tau\left(  a_{1}^{\frac{-1}{2}}R_{1,2,1,2}a_{2}^{\frac{-1}{2}}\right)  =0,
\]
either when each $a_{j}$ commutes with its first-order derivatives
$\partial_{1}a_{j}$ and $\partial_{2}a_{j}$ (but $a_{1}$ and $a_{2}$ need not
commute in general) or when $a_{1}a_{2}=1$, where $\tau$ is the unique
faithful tracial state of $A_{\theta}$.

There are many concrete nontrivial cases satisfying the criteria of this
theorem. For example, each of $a_{j}:=u_{j}+3+u_{j}^{\ast}$ for $j=1,2$
commutes with $\partial_{k}a_{j}=2\pi i\delta_{jk}\left(  u_{j}+u_{j}^{\ast
}\right)  $ for $k=1,2$ (but $a_{1}$ and $a_{2}$ do not commute). More
generally, one can take $a_{j}:=f\left(  u_{j}\right)  $ for a holomorphic
function $f$ on a neighborhood of $\mathbb{T}$ with $f|_{\mathbb{T}}$ positive-valued.

After completing this paper, we became aware of an interesting paper by D\c
{a}browski and Sitarz \cite{DaSi} that derives a Gauss-Bonnet theorem in the
framework of spectral triples and Dirac operators for a class of non-conformal
deformations of the flat two-torus. Also, it should be noted that recent work
of Arnlind and Wilson \cite{ArWi,ArWi2} discusses noncommutative Levi-Civita
connections in a framework with certain similarities to Rosenberg's.

%

{\medskip}%

\noindent Mira A. Peterka \newline\noindent Instytut Matematyczny, Polska
Akademia Nauk, ul. Sniadeckich 8, Warszawa, 00-956 Poland \newline\noindent peterka@math.berkeley.edu%

{\medskip}%

\noindent Albert J. L. Sheu \newline\noindent Department of Mathematics,
University of Kansas, Lawrence, KS 66045, U. S. A. \newline\noindent asheu@ku.edu

\end{document}